\documentclass[14pt]{extreport}
\usepackage[cp1251]{inputenc}
\usepackage[ukrainian,russian,english]{babel}
\usepackage{amsmath,amsfonts,amssymb}
\usepackage{latexsym,hhline,graphics,epsfig}

\newcounter{theorem}

\newcounter{lemma}
\newcommand{\be}{\begin{equation}}
\newcommand{\ee}{\end{equation}}

\begin{document} \begin{flushleft}

\bigskip

\begin{center}{\large\bf Galina Bakhtina, Inna Dvorak, Iryna Denega}\end{center}

\begin{center}{\large\bf Problems on extremal decomposition of the complex plane}\end{center}
\end{flushleft}


\begin{flushright}
\parbox{16.5cm}{\footnotesize Considered here is one quite general problem about description 
of extremal configurations maximizing the product of
inner radii non-overlapping domains.}
\end{flushright}

\bigskip

Let $\mathbb{N}$, $\mathbb{R}$ be a set of natural and real numbers,
respectively, $\mathbb{C}$ be a complex plane,
$\overline{\mathbb{C}}=\mathbb{C}\bigcup\{\infty\}$ be a
Riemann sphere and $\mathbb{R^{+}}=(0,\infty)$.

Let $r(B,a)$ be a inner radius of domain
$B\subset\overline{\mathbb{C}}$, with respect to a point $a\in B$
(see \cite{1,2,heim,Dubinin-1994,Bakh-Bakhtina-Zel-2008}).

Let $n\in \mathbb{N}$. A set of points $A_{n}:=\left\{a_{k} \in
\mathbb{C}:\,k=\overline{1,n}\right\},$ is called \emph{$n$-radial}
system, if $|a_{k}|\in\mathbb{R^{+}}$, $k=\overline{1,n}$, and
$0=\arg a_{1}<\arg a_{2}<\ldots<\arg a_{n}<2\pi$.

Denote $$\alpha_{k}:=\frac{1}{\pi}\arg \frac{a_{k+1}}{a_{k}},\,
\alpha_{n+1}:=\alpha_{1},\,k=\overline{1,n}.$$

Consider next functional
\begin{equation}\label{2funk}J_n(\gamma)=r^\gamma\left(B_0,0\right)
\prod\limits_{k=1}^n r\left(B_k,a_k\right),\end{equation} where
$\gamma\in\mathbb{R^{+}}$, $B_{0}$, $B_{1}$, $B_{2}$,..., $B_{n}$,
$n\geq2$ is non-overlapping domains (e.i. $B_{p}\cap B_{j}={\O}$ if
$p\neq j$) in $\overline{\mathbb{C}}$, $a_{0}=0$, $|a_{k}|=1$,
$k=\overline{1,n}$, $a_{k}\in B_{k}$, $k=\overline{0, n}$ and
$\gamma\leq n$.

One of the main concepts in this paper is a notion of quadratic
differential. Quadratic differential is a convenient tool of
describing extremals on the plane. A great collection of definitions and results about quadratic
differentials one can find in monographs by well known analyst
J. Jenkins \cite{Djen-1962}.

\begin{figure}[h]
\begin{center}{\includegraphics[width=6cm]{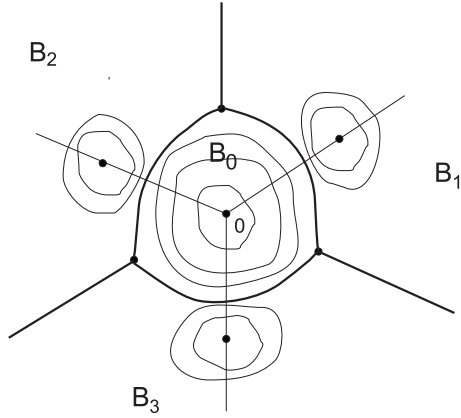}}\end{center}
\caption{Structure of the trajectories of the quadratic
differential (\ref{2kv}) for $n=3$ and $\gamma=1$.}\label{void23}
\end{figure}

For example, we consider next quadratic
differential
\begin{equation}\label{2kv}Q(w)dw^2=-\frac{(n^2-\gamma)w^n+\gamma}{w^2(w^n-1)^2}.\end{equation}
If $n=3$ and
$\gamma=1$ then (\ref{2kv}) has four circular domains $B_{0}$, $B_{1}$, $B_{2}$, $B_{3}$ such that
$\overline{B_{0}\cup B_{1}\cup B_{2}\cup B_{3}}=\overline{\mathbb{C}}$. 
A simple poles of the quadratic differential (\ref{2kv}) are
$w_{1}=-\frac{1}{2}$,
$w_{2}=\frac{1}{4}+\frac{\sqrt{3}}{4}i$,
$w_{3}=\frac{1}{4}-\frac{\sqrt{3}}{4}i$. The structure of the trajectories of (\ref{2kv})
is shown in Fig. \ref{void23}.

\newpage

Let $$F_{\delta}(x)=2^{x^2+6}\cdot x^{x^2+2-2\delta}\cdot
(2-x)^{-\frac{1}{2}(2-x)^2} (2+x)^{-\frac{1}{2}(2+x)^2},$$ $$x \in
(0,2],\quad 0\leqslant\delta\leqslant0,7.$$

We proved 

\textbf{Theorem.} \emph{Let $n\in \mathbb{N}$, $n\geqslant 5$,
$\gamma\in (0; 1,75\,]$, $0\leqslant\delta\leqslant0,7$. Then for any
$n$-radial system of points $A_n=\{a_k\}_{k=1}^n$, $|a_k|=1$, and
any system of non-overlapping domains $B_k$, $a_k\in
B_k\subset\overline{\mathbb{C}}$, $k=\overline{0,n}$, $a_{0}=0,$
we have inequality
$$J_n(\gamma)\leqslant\gamma^{-\frac{\delta\cdot n}{2}}\cdot\left(
\prod^n_{k=1} \alpha_k\right)^{\delta} \cdot \left[F_{\delta}
\left(\frac{2}{n}\sqrt{\gamma}\right)\right]^{\frac{n}{2}},$$
where equality holds if $a_k$ and $B_k$, $k=\overline{0,n}$, are, respectively, poles and circular domains of
the quadratic differential} (\ref{2kv}).

Method of proof of this theorem
is based on using of separating transformation \cite{Dubinin-1988,Dubinin-1994} and ideas of works
\cite{kovalev,Bakh-Bakhtina-Zel-2008,de1,de2}.

In an analogous way as in papers \cite{de1,de2,Bahtin} we obtain next inequality for functional (\ref{2funk})
\begin{equation}\label{122a}J_n(\gamma)\leqslant\gamma^{-\frac{\delta n}{2}}\left(
\prod\limits^n_{k=1} \alpha_k\right)^{\delta}
\left[\prod\limits^n_{k=1}
F_{\delta}\left(\alpha_{k}\sqrt{\gamma}\right)\right]^{1/2},\,\delta\in[0;
0,7].\end{equation} Consider the functional
$$\tilde J_n(\gamma)=\gamma^{\frac{\delta n}{2}}\left(\prod\limits^n_{k=1}
\alpha_k\right)^{-\delta}J_n^{(\gamma)}.$$ It follows from (\ref{122a}) that
$$\tilde J_n(\gamma)\leqslant\left[\prod\limits^n_{k=1}
F_{\delta}\left(\alpha_{k}\sqrt{\gamma}\right)\right]
^{\frac{1}{2}}.$$
Further, consider an extremal problem $$\prod^n_{k=1}F_{\delta}(x_k)
\longrightarrow \max, \quad \sum^n_{k=1}x_k=2\sqrt{\gamma}, \quad x_k=\alpha_{k}\sqrt{\gamma}$$
$$0<x_k\leqslant 2, \quad 0\leqslant\delta\leqslant 0,7.$$
Let $\Psi_{\delta}(x)=\ln \left(F_{\delta}(x)\right)$. And let
$X^{(0)}=\left\{x_k^{(0)}\right\}_{k=1}^{n}$ be an arbitrary
extremal point of this problem.

\begin{figure}[h]
\begin{center}{\includegraphics[width=6cm]{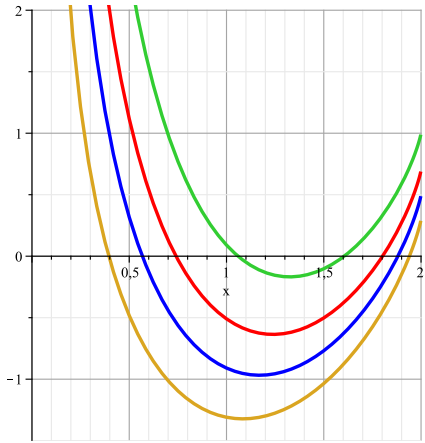}}\end{center}
\caption{Graph of the function $\Psi_{\delta}'(x)$}\label{ris2}
\end{figure}

It follows from the paper \cite{kovalev} that:
if $0<x_k^{(0)}<x_j^{(0)}<2$, $k\neq j$, than
$$\Psi'_{\delta}(x_k^{(0)})=\Psi'_{\delta}(x_j^{(0)}),$$
where $k,j=\overline{1,n}$, $k\neq j$, $0\leqslant\delta\leqslant0,7$,
$$\Psi_{\delta}'(x)=2x\ln(2x)+(2-x)\ln(2-x)
-(2+x)\ln(2+x)+\frac{2}{x}-\frac{2\delta}{x}$$
(see Fig. \ref{ris2}).

We will show that the following assertion is true
$$x_{1}^{(0)}=x_{2}^{(0)}=\ldots=x_{n}^{(0)}.$$
Let
$$\sigma_1:=\sigma_1(\delta,\gamma)=\min\limits_{1\leqslant k\leqslant n}x_{k}^{(0)}(\delta,\gamma),$$
$$\sigma_0:=\sigma_0(\delta,\gamma)=\max\limits_{1\leqslant k\leqslant n}x_{k}^{(0)}(\delta,\gamma),$$
$\sigma_1\leqslant\sigma_k\leqslant\sigma_0$, $k=\overline{1,n}$,
$0\leqslant\delta\leqslant 0,7$, $\gamma\in(0; 1,75]$.

A function $\Psi_{\delta}''(x)$ is strictly increasing on $(0,2)$ for each fixed $\delta$. Thus
$$sign \Psi_{\delta}''(x)\equiv sign (x-x_{0}(\delta,\gamma)).$$

It is not difficult to see if $\sigma_0\leqslant x_{0}(\delta,\gamma)$ then from the conditions of the problem we obtain $x_{1}^{(0)}=\ldots=x_{n}^{(0)}$.

Assume $x_{0}(\delta,\gamma)\leqslant\sigma_0<1,62$. Then
$$\sigma_1\leqslant\frac{1}{n-1}\sum^{n-1}_{k=1}x_k^{(0)}=
\frac{2\sqrt{\gamma}-\sigma_0}{n-1}\leqslant\frac{2-\sigma_{0}}{n-1}.$$
Therefore for $n\geqslant 5$, $\gamma\in(0; 1,75]$,
$0\leqslant\delta\leqslant 0,7$, inequality holds
$$\sigma_{1}\leqslant(2\sqrt{\gamma}-\sigma_{0})/(n-1) \leqslant(2,645751-x_{0}(\delta,\gamma))/4\leqslant$$$$\leqslant(2,645751-1,08441)/4<0,390335.$$
Thus
$$\Psi_{\delta}'(\sigma_1)>\Psi_{\delta}'(0,390335)>\Psi_{0,7}'(0,390335)=0,027069>$$$$>0,018707=
\Psi_{0}'(1,62)\geqslant \Psi_{\delta}'(\sigma_0).$$ Hence $\sigma_0 \notin [x_{0}(\delta,\gamma);\,1,62)$.

Let $1,62\leqslant\sigma_0\leqslant2$. Then for $n\geqslant5$, $\gamma\in(0; 1,75]$,
$0\leqslant\delta\leqslant 0,7$, we have
$$\sigma_{1}\leqslant(2\sqrt{\gamma}-\sigma_{0})/(n-1)\leqslant(2,645751-1,62)/4<0,256438.$$ That is
$$\Psi_{\delta}'(\sigma_1)>\Psi_{\delta}'(0,256438)>\Psi_{0,7}'(0,256438)=1,130326>$$$$>1=
\Psi_{0}'(2)\geqslant \Psi_{\delta}'(\sigma_0).$$ So, $\sigma_0 \notin [x_{0}(\delta,\gamma);\,2)$.

Therefore, extremal set $\left\{x_{k}^{(0)}\right\}_{k=1}^{n}$
is possible only if $$x_{1}^{(0)}=x_{2}^{(0)}=\ldots=x_{n}^{(0)}.$$
From the foregoing the following relation holds $$\prod\limits^n_{k=1}
F_{\delta}\left(\alpha_{k}\sqrt{\gamma}\right)\leqslant\left[F_{\delta}
\left(\frac{2}{n}\sqrt{\gamma}\right)\right]^{n}.$$ A sign of equality is verified directly.
Theorem is thus proved.

{\small
\renewcommand{\refname}{References}

\end{document}